\documentclass [11pt] {article}

\usepackage{amsmath,amssymb}

\newtheorem{theorem}{Theorem}

  \newtheorem{corollary}{Corollary}
  
  \newtheorem{remark}{Remark}

\begin{document}

\title{{\Large{\bf On Hall's conjecture}}}

\author{Andrej Dujella}

\date{}
\maketitle

\begin{abstract}
We show that for any even positive integer $\delta$ there exist polynomials
$x$ and $y$ with integer coefficients such that
$\deg(x)=2\delta$, $\deg(y)=3\delta$ and $\deg(x^3-y^2)=\delta+5$.

\end{abstract}

\footnotetext{


The author was supported
by the Ministry of Science, Education and Sports, Republic of Croatia, grant
037-0372781-2821.
}


\bigskip

Hall's conjecture asserts that for any $\varepsilon >0$, there exists a
constant $c(\varepsilon)>0$ such that if $x$ and $y$ are positive integers satisfying
$x^3-y^2 \neq 0$, then $|x^3-y^2| > c(\varepsilon) {x}^{1/2-\varepsilon}$.
It is known that Hall's conjecture follows from the $abc$-conjecture. For a stronger
version of Hall's conjecture which is equivalent to the $abc$-conjecture see \cite[Ch. 12.5]{BG}.
Originally, Hall \cite{Hall} conjectured that there is $C>0$ such that $|x^3-y^2| \geq C \sqrt{x}$
for positive integers $x,y$ with $x^3-y^2\neq 0$, but this formulation is unlikely to be true.
Danilov \cite{Dan} proved that $0<|x^3-y^2| < 0.97 \sqrt{x}$ has infinitely many solutions
in positive integers $x,y$; here $0.97$ comes from $54\sqrt{5}/125$. For examples
with ``very small'' quotients $|x^3-y^2|/\sqrt{x}$, up to $0.021$, see \cite{Elk} and \cite{JCHS}.

It is well known that for non-constant complex polynomials $x$ and $y$,
such that $x^3\neq y^2$, we have
$\deg(x^3 - y^2) / \deg(x) > 1/2$.
More precisely, Davenport \cite{Dav} proved that for such polynomials the inequality
\begin{equation} \label{dav}
\deg(x^3-y^2) \geq \tfrac{1}{2} \deg(x) + 1
\end{equation}
holds. This statement also follows from Stothers-Mason's $abc$ theorem for polynomials
(see, e.g., \cite[Ch. 4.7]{La}).
Zannier \cite{Zan2} proved that for any positive integer $\delta$ there exist complex polynomials
$x$ and $y$ such that $\deg(x)=2\delta$, $\deg(y)=3\delta$ and
$x,y$ satisfy the equality in Davenport's bound (\ref{dav}). In his previous paper \cite{Zan1},
he related the existence of such examples with coverings of the Riemann sphere,
unramified except above $0$, $1$ and $\infty$.

It is natural to ask whether examples with the equality in (\ref{dav}) exist for
polynomials with integer (rational) coefficients. Such examples are known only
for $\delta=1,2,3,4,5$ (see \cite{B-S,Elk}).
The first example for $\delta=5$ was found by Birch, Chowla, Hall and Schinzel \cite{BCHS}.
It is given by
$$ x= \frac{t}{9}(t^9+6t^6+15t^3+12), \quad y=\frac{1}{54}(2t^{15}+18t^{12}+72t^{9}+144t^{6}+135t^{3}+27), $$
while then
$$ x^3-y^2 = -\frac{1}{108}(3t^{6}+14t^{3}+27) $$
(note that $x,y$ are integers for $t\equiv 3\pmod{6}$).
One more example for $\delta=5$ has been found by Elkies \cite{Elk}:
$$ x = t^{10} - 2t^9 +33t^8 - 12t^7 + 378t^6 + 336t^5 + 2862t^4 +2652t^3 + 14397t^2 +9922t+ 18553,
\vspace{-2ex} $$
{\small
\begin{eqnarray*}
y &\!\!=\!\!& t^{15} - 3t^{14} + 51t^{13} -67t^{12} + 969t^{11} +33t^{10}
+10963t^9 + 9729t^8 + 96507t^7 \\
& &\mbox{}+ 108631t^6
+ 580785t^5 + 700503t^4 + 2102099t^3 + 1877667t^2 + 3904161t+ 1164691, \vspace{-5ex}
\end{eqnarray*} \vspace{-5ex}
\begin{eqnarray*}
 x^3-y^2 &\!\!=\!\!& 4591650240t^6-5509980288t^5+101934635328t^4+58773123072t^3\\
 & & \mbox{}+730072388160t^2+1151585880192t+5029693672896.
\end{eqnarray*} }%
In these examples we have
$$ \deg(x^3 - y^2) / \deg(x) = 0.6, $$
and it seems that no examples of polynomials with integer coefficients,
satisfying $x^3-y^2\neq 0$ and $ \deg(x^3 - y^2) / \deg(x) < 0.6$,
were published until now.

In this note we will show the following result.

\begin{theorem} \label{tm1}
For any $\varepsilon>0$ there exist
polynomials $x$ and $y$ with integer coefficients
such that $x^3\neq y^2$ and
$\deg(x^3 - y^2) / \deg(x) < 1/2 + \varepsilon$.

More precisely, for any even positive integer $\delta$ there exist polynomials
$x$ and $y$ with integer coefficients such that
$\deg(x)=2\delta$, $\deg(y)=3\delta$ and $\deg(x^3-y^2)=\delta+5$.
\end{theorem}

As an immediate corollary we obtain a nontrivial lower bound for the number of integer solutions
to the inequality $|x^3-y^2| < x^{1/2+\varepsilon}$ with $1\leq x \leq N$ (heuristically, it is
expected that this number is around $N^{\varepsilon}$).

\begin{corollary} \label{cor1}
For any $\varepsilon >0$ and positive integer $N$ by
$\mathcal{S}(\varepsilon,N)$ we denote the number of integers $x$, $1\leq x\leq N$,
for which there exists an integer $y$ such that $0 < |x^3-y^2| < x^{1/2+\varepsilon}$.
Then we have
$$ \mathcal{S}(\varepsilon,N) \gg N^{\varepsilon / (5+4\varepsilon)}. $$
\end{corollary}

Indeed, take $\delta$ to be the smallest even integer greater that $5/(2\varepsilon)$,
so that $5/(2\varepsilon) < \delta < 5/(2\varepsilon)+2$, and take $x=x(t)$, $y=y(t)$
as in Theorem \ref{tm1}. Then for sufficiently large $t$ we have $x=O(t^{2\delta})$ and
$ |x^3-y^2| = O(t^{\delta +5})=O(x^{\frac{1}{2}+\frac{5}{2\delta}}) < x^{1/2+\varepsilon}.$
Therefore, \vspace{-1.0ex}
$$ \mathcal{S}(\varepsilon,N) \gg N^{1/(2\delta)}
\gg N^{\varepsilon/(5+4\varepsilon)}. $$

\medskip

Here is an explicit example which improves the quotient $\deg(x^3 - y^2) / \deg(x) = 0.6$
from the above mentioned examples
by Birch, Chowla, Hall, Schinzel and Elkies, as
$\deg(x^3 - y^2) / \deg(x) = 31/52 = 0.5961...$ : \vspace{-2.3ex}

{\footnotesize
\begin{eqnarray*}
& x & \hspace*{-0.2cm} = 281474976710656t^{52}+3799912185593856t^{50}+24189255811072000t^{48}+96537120918732800t^{46}\\
 & &\mbox{}\hspace*{-0.4cm}+270892177293312000t^{44}+568175382432317440t^{42}+924393098014883840t^{40}\\
 & &\mbox{}\hspace*{-0.4cm}+1194971570896896000t^{38}+1247222961904025600t^{36}+1062249296822272000t^{34}\\
 & &\mbox{}\hspace*{-0.4cm}+743181990714408960t^{32}+428630517911388160t^{30}+203971125837824000t^{28}+100663296t^{27}\\
 & &\mbox{}\hspace*{-0.4cm}+79960271015116800t^{26}+729808896t^{25}+25720746147840000t^{24}+2359296000t^{23}\\
 & &\mbox{}\hspace*{-0.4cm}+6745085391667200t^{22}+4482662400t^{21}+1428736897843200t^{20}+5554176000t^{19}\\
 & &\mbox{}\hspace*{-0.4cm}+241375027200000t^{18}+4706795520t^{17}+31982191104000t^{16}+2782494720t^{15}+3250264320000t^{14}\\
 & &\mbox{}\hspace*{-0.4cm}+1148928000t^{13}+245895686400t^{12}+326476800t^{11}+13292822400t^{10}+61776000t^{9}+484380000t^{8}\\
 & &\mbox{}\hspace*{-0.4cm}+7344480t^{7}+10894000t^{6}+496080t^{5}+130625t^{4}+15750t^{3}+629t^{2}+150t+4, \\
\end{eqnarray*} \vspace{-7.0ex}
\begin{eqnarray*}
& y & \hspace*{-0.2cm} = 4722366482869645213696t^{78}+95627921278110315577344t^{76}+931486788746037518401536t^{74}\\
 & &\mbox{}\hspace*{-0.4cm}+5812273909720700361375744t^{72}+26102714713365300532740096t^{70}+89873242715073754863501312t^{68}\\
 & &\mbox{}\hspace*{-0.4cm}+246761827996223603178733568t^{66}+554869751478978106456276992t^{64}\\
 & &\mbox{}\hspace*{-0.4cm}+1041377162422256031202541568t^{62}+1654256777803799676753805312t^{60}\\
 & &\mbox{}\hspace*{-0.4cm}+2247766244734980591395536896t^{58}+2633529391786763986554322944t^{56}\\
 & &\mbox{}\hspace*{-0.4cm}+2676840149412734907329806336t^{54}+2533274790395904t^{53}+2371433108159248512627769344t^{52}\\
 & &\mbox{}\hspace*{-0.4cm}+35465847065542656t^{51}+1837294956807449113993936896t^{50}+234486247786020864t^{49}\\
 & &\mbox{}\hspace*{-0.4cm}+1247823926411289395000770560t^{48}+973569167884025856t^{47}+743994544482135039635619840t^{46}\\
 & &\mbox{}\hspace*{-0.4cm}+2847272221544546304t^{45}+389682593956278112836648960t^{44}+6236328797675716608t^{43}\\
 & &\mbox{}\hspace*{-0.4cm}+179279686440609529032867840t^{42}+10618254681610125312t^{41}+72388134028773255869890560t^{40}\\
 & &\mbox{}\hspace*{-0.4cm}+14399046085119049728t^{39}+25611943886548098204303360t^{38}+15806610071787405312t^{37}\\
 & &\mbox{}\hspace*{-0.4cm}+7922395450159324505047040t^{36}+14200560742834372608t^{35}+2135839807968003238133760t^{34}\\
 & &\mbox{}\hspace*{-0.4cm}+10514148446410113024t^{33}+499883693495498613719040t^{32}+6441026076788391936t^{31}\\
 & &\mbox{}\hspace*{-0.4cm}+101073262762096181903360t^{30}+3269189665642512384t^{29}+17550157782838363029504t^{28}\\
 & &\mbox{}\hspace*{-0.4cm}+1373442845007937536t^{27}+2598168579136061177856t^{26}+476068223096193024t^{25}\\
 & &\mbox{}\hspace*{-0.4cm}+325093317533140516864t^{24}+135395930768670720t^{23}+34019036843474681856t^{22}\\
 & &\mbox{}\hspace*{-0.4cm}+31339645700014080t^{21}+2939255644452962304t^{20}+5838612910571520t^{19}+206402445920944128t^{18}\\
 & &\mbox{}\hspace*{-0.4cm}+862650209710080t^{17}+11551766627438592t^{16}+99129281310720t^{15}+502656091170048t^{14}\\
 & &\mbox{}\hspace*{-0.4cm}+8633278321920t^{13}+16468534726592t^{12}+550276346880t^{11}+389483950128t^{10}+24450210720t^{9}\\
 & &\mbox{}\hspace*{-0.4cm}+6312333144t^{8}+705350880t^{7}+68685241t^{6}+11812545t^{5}+642429t^{4}+94050t^{3}+6591t^{2}+225t+19, \\
\end{eqnarray*}
\begin{eqnarray*}
& x^3-y^2 & \hspace*{-0.2cm} = -905969664t^{31}-8380219392t^{29}-35276193792t^{27}-89379569664t^{25}-151909171200t^{23}\\
 & &\mbox{}\hspace*{-0.8cm}-182680289280t^{21}-159752355840t^{19}-102786416640t^{17}-48661447680t^{15}-16772918400t^{13}\\
 & &\mbox{}\hspace*{-0.8cm}-4116359520t^{11}-692649360t^{9}-75171510t^{7}-297t^{6}-4749570t^{5}-891t^{4}-144450t^{3}-891t^{2}\\
 & &\mbox{}\hspace*{-0.8cm}-1350t-297.
\end{eqnarray*} }%

Now we describe the general construction. Let us define the binary recursive sequence by
$$ a_1=0, \quad a_2=t^2+1, \quad a_m = 2ta_{m-1}+a_{m-2}.$$
Thus, for $m\geq 2$, $a_m$ is a polynomial
in variable $t$, of degree $m$.
Put $u = a_{k-1}$ and $v = a_{k}$ for an odd positive integer $k\geq 3$.
We search for examples with $x=O(v^2)$, $y=O(v^3)$ and $x^3-y^2=O(v)$.
Note that
\begin{equation} \label{dioph}
v^2-2tuv-u^2 = -(a_2^2-2ta_1a_2-a_1^2)= -(t^2+1)^2.
\end{equation}
Therefore, we may take
\begin{eqnarray*}
x &\!\!=\!\!& av^2+buv+cu+dv+e, \\
y &\!\!=\!\!& fv^3+gv^2u+hv^2+iuv+ju+mv+n,
\end{eqnarray*}
with unknown coefficients $a,b,c,\ldots,n$, which will be
determined so that in the expression for $x^3-y^2$
the coefficients with $v^6, uv^5, v^5, \ldots , v^2, uv$
are equal to $0$.
We find the following (polynomial) solution:
\begin{eqnarray*}
x &\!\!=\!\!& v^2-2tuv+6v-6tu+(t^4+5t^2+4), \\
y &\!\!=\!\!& -2tv^3+(4t^2+1)uv^2 -9tv^2+(18t^2+9)uv +(-2t^5-4t^3-2t)v \\
& &\mbox{}+(t^4+20t^2+19)u+ (-9t^5-18t^3-9t).
\end{eqnarray*}
Using (\ref{dioph}), it is easy to check that we have
$$ x^3-y^2 = -27(t^2+1)^2(2v-2tu+11t^2+11). $$
Therefore, $\deg(x)=2k-2$ and $\deg(x^3-y^2)=k+4$. Also,
$$ \deg(x^3-y^2) / \deg(x) = (k+4)/(2k-2), $$
which tends to $1/2$ when $k$ tends to infinity.
The above explicit example corresponds to $k=27$.

\bigskip

Comparing with Davenport's bound, our polynomial $x$ and $y$ satisfy
$$ \deg(x^3-y^2) = \tfrac{1}{2} \deg(x) + 5. $$
Thus, although our examples $(x,y)$ do not give the equality in Davenport's bound (\ref{dav}),
they are very close to the best possible result for $\deg(x^3-y^2)$, and it seems
that this is the first known result of the form that $\deg(x^3-y^2) - \tfrac{1}{2} \deg(x)$
is bounded by an absolute constant, for polynomials $x,y$ with integer coefficients
and arbitrarily large degrees.

Since $(t^2+1)$ divides $a_m$ for all $m$, it could be noted that $(t^2+1)$ divides $x$ and $(t^2+1)^2$
divides $y$. Hence, with $x=(t^2+1)X$ and $y=(t^2+1)^2Y$, we have
$$ \deg(X^3 - (t^2+1)Y^2) = \tfrac{1}{2} \deg(X). $$
This shows that the only branch points of the rational function $x^3/y^2$ are $0$, $1$ and $\infty$,
which is in agreement with the results of Zannier \cite{Zan1,Zan2}.

\bigskip

Let us give an interpretation of our result in terms of polynomial Pell's equations.
Following a suggestion by N.~Elkies, we put $v-tu = (t^2+1)z$.
Then the expressions of $x$ and $x^3-y^2$ simplify considerably, and
we get
$x = (t^2+1)(z^2+6z+4)$, $x^3 - y^2 = -27(t^2+1)^3 (2z+11)$ which gives
$y^2=(t^2+1)^3(z^2+1)(z^2+9z+19)^2$. Thus, we need that $z^2+1= (t^2+1)w^2$, i.e
\begin{equation} \label{pel}
z^2 - (t^2+1)w^2 = -1.
\end{equation}
The fundamental solution of Pell's equation (\ref{pel}) is $(z,w)=(t,1)$. Taking $t=z$, we obtain the
identity
$$ (z^2+6z+4)^3 - (z^2+1) (z^2+9z+19)^2 = - 27(2z+11), $$
which is equivalent to Danilov's example \cite{Dan} (and by taking $z^2+1=5w^2$
and $2z+11 \equiv 0 \pmod{125}$,
we get a well-known sequence of numerical examples with $|x^3-y^2| < \sqrt{x}$).

However, if we consider (\ref{pel}) as a polynomial Pell's equation (in variable $t$), we obtain
the sequence of solutions
$$ z_1=t, \quad z_2=4t^3+3t, \quad z_{k}=(4t^2+2)z_{k-1}-z_{k-2} . $$
This gives exactly the sequences of polynomials $x$ and $y$,
as given above.

\bigskip

\begin{remark}
{\rm In \cite{Dan2}, Danilov consireded small values of $|x^4-Ay^2|$, for integers $A$ satisfying certain
conditions. Using the formula
\begin{equation} \label{eq:Dan2}
 (27z+7)^4 - (81z+20)^2 \cdot \frac{(81z+22)^2+2}{81} = 4z+1,
\end{equation}
he proved that if the Pellian equation $u^2-81Av^2=-2$ has a solution, then the inequality
$|x^4-Ay^2| < \frac{4}{27}|x|$ has infinitely many integer solutions $x,y$.
By applying a similar construction, as above, to
Danilov's formula (\ref{eq:Dan2}), we obtain
the sequences $x_k$ and $y_k$ of polynomials in variable $t$ with $\deg(x_k)=2k+1$, $\deg(y_k)=4k$ and
$\deg(x^4 - (t^2+2)y^2) = \deg(x)=2k+1$.
For example, for $k=3$ we have
{\small
\begin{eqnarray*}
 x &\!\!=\!\!& 8t^7+28t^5+28t^3+7t-1, \\
 y &\!\!=\!\!& 64t^{13}+384t^{11}+880t^{9}+960t^{7}-16t^{6}+504t^{5}-40t^{4}+112t^{3}-24t^{2}+7t-2,
\end{eqnarray*} }%
and then
$$ x^4-(t^2+2)y^2= 32t^7+112t^5+112t^3+28t-7. $$ }
\end{remark}

\bigskip

{\bf Acknowledgements.} The author is grateful to Yann Bugeaud, Noam
Elkies, Clemens Fuchs, Boris \v{S}irola and Umberto Zannier for their very interesting and useful
comments on the previous version of this note.

\medskip

{\footnotesize
Department of Mathematics, University of Zagreb,

 Bijeni\v{c}ka cesta 30, 10000 Zagreb, Croatia

{\em E-mail address}: {\tt duje@math.hr} }

\end{document}